\newcommand{\me}{\mathrm{e}}
\newcommand{\mi}{\mathrm{i}}

\documentclass[12pt,a4paper]{article}
\usepackage{empheq}
\usepackage{amsmath}
\usepackage{amssymb}
\usepackage{tikz}
\usepackage{amsfonts}
\usepackage{amsthm}
\usepackage{xcolor}
\theoremstyle{plain}
\newtheorem*{theorem*}{Theorem}
\begin{document}

\title{Far-Field Sensitivity to Local Boundary Perturbations in 2D Wave Scattering}

\author{Erik~Garc\'ia~Neefjes${}^{1}$\footnote{Author for correspondence. E-mail: \texttt{erik.garcia@mq.edu.au}},   Stuart~C.~Hawkins${}^{1}$
\\[5pt]
{\footnotesize
${}^{1}$ School of Mathematical and Physical Sciences, Macquarie University, Sydney NSW 2109, Australia}
}

% Use default \verb|\maketitle|.
\date{\today}
\maketitle

\begin{abstract}
We numerically investigate the sensitivity of the scattered wave field to perturbations in the shape
of a scattering body illuminated by an incident plane wave.
This study is motivated by recent work on the inverse problem of reconstructing a scatterer
shape from measurements of the scattered wave at large distances from the scatterer.
For this purpose we consider
star-shaped scatterers represented using cubic splines, and our approach is based on
a Nystr\"om method-based discretisation of the
shape derivative.
Using the singular value decomposition,
we identify fundamental geometric modes that most strongly influence the scattered wave, providing insight into the most visible boundary features in scattering data.
\end{abstract}

% By default we include a table of contents in each paper.
\tableofcontents

% Use \verb|\section|, \verb|\subsection|, \verb|\subsubsection| and 
% possibly \verb|\paragraph| to structure your document.  Make sure 
% you \verb|\label| them for cross-referencing with \verb|\ref|.
\section{Introduction}
\label{sec:intro}

When a wave interacts with a \textit{scatterer}, its properties are modified through a process known as \textit{wave scattering} \cite{born2013principles}. The {scattered wave} pattern depends strongly on the shape and properties of the scatterer and the direction from which the wave approaches. Of particular interest is the response of the scattered wave at large distances from the scatterer, commonly known as the \textit{far-field}. In this region, the scattered wave approximates an outgoing spherical wave, with its amplitude decreasing inversely with distance from the scatterer. Far-field patterns are important in practice since often measurements must be made at considerable distances from the scattering object. Given far-field measurements, reconstructing properties of the scatterer is a classical \textit{ill-posed} inverse problem. While, consequently, the forward map from scatterers to far-fields has no inverse, various methods including iterative, qualitative, and regularization approaches have been developed to obtain useful approximate solutions \cite{colton2019inverse}. 

For the inverse shape problem, remarkable reconstructions have been obtained with (deterministic) regularized Newton-type methods \cite{borges2022multifrequency}, as well as (stochastic) Bayesian counterparts \cite{yang2021bayesian}. {In tackling this problem, following our work in \cite{ganesh2020efficient}, we noticed some interesting phenomena related to how well certain kinds of scatterer-dependent features can be reconstructed. In this paper, we aim to investigate these phenomena further by considering the forward problem for an illustrative geometry, and show how our previous observations can be explained via analysis of the shape-derivative. Interpreting this derivative is facilitated by Singular Value Decomposition (SVD), which allows identification of the geometric perturbation modes that are most \textit{visible} in the far-field, and hence expected to be \textit{reconstructible} in the inverse problem configuration.}

Unlike traditional approaches which use global basis functions like \textit{Karhunen-Lo\`eve} expansions \cite{yang2021bayesian, ganesh2020efficient}, the use of cubic splines allows for localised geometric perturbations, which is particularly valuable for understanding how specific geometric features influence the far-field. However, for such analysis it is paramount to be able to generate accurate numerical solutions, especially for complicated shapes with high spatial resolution. We overcome this by using the spectral \textit{Nystr\"om} method \cite{colton2019inverse} for solving Boundary Integral Equations (BIEs), and we ensure the required smoothness of our scatterer's boundary representation.

In Section \ref{Section:Scatterer representation} we present the the class of scatterers considered in this article, Section \ref{section: Model} focuses on the wave scattering model including the governing BIEs that lead to the sensitivity matrix. In Section \ref{section:results} we demonstrate our corresponding numerical results with focus on the SVD analysis. Conclusions are given in Section \ref{section:conclusions}.
\begin{figure}
\centering
\includegraphics[scale=0.15]{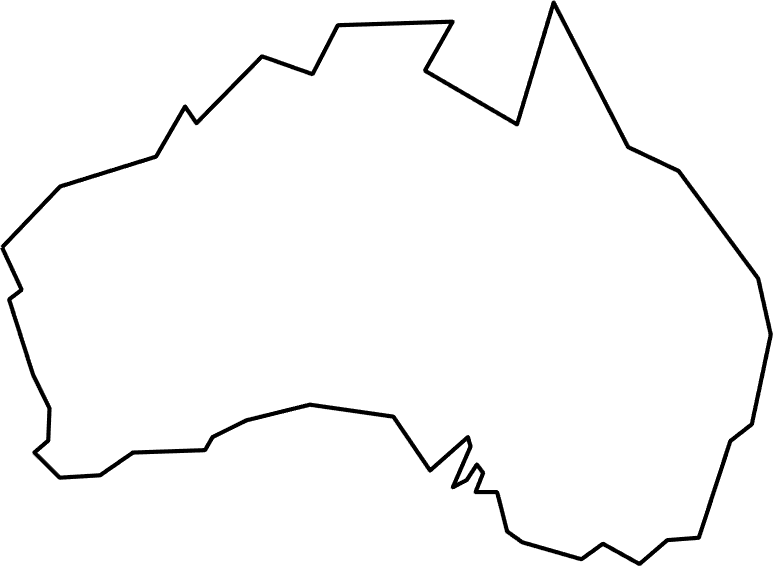}
\caption{Visualisation of our reference scatterer, a polygonal approximation to the coast
of mainland Australia.}
\label{fig:Aus True Cartoon}
\end{figure}

%The far-field also exhibits simpler mathematical properties than the complete scattered field, making it a valuable tool for analyzing scattering problems, though this simplification comes with inherent limitations in the information it can provide about the scatterer.

% The Introduction has to show that your story is worth telling in
% detail.  The Introduction is likely to be all an interested reader
% reads, again it must be complete in itself.  Place your work in the
% context of other research.  Summarise your main results, albeit in a
% suitably simplified form.

% Face it: only the dedicated are going to want to wade through the 
% details of the rest of the paper.  Give the key points in your 
% Introduction.

%\section{Write well}
%\label{sec:ww}

%Be definite.  Be descriptive.  Be precise.  Cross reference.  Use
%short sentences.

% Keep close together nouns and their verbs: that is, write ``the cat
% sat on the mat'' not ``the cat on the mat sat''.

% Structure your writing using the \emph{rule of three}.  Each paragraph 
% is to make a point: overview the paragraph in the first sentence; 
% develop the point in the body; and summarise the paragraph in the last 
% sentence.  Each section develops a theme: overview a section in its 
% first paragraph; develop the theme in the body; and summarise in its 
% last paragraph.  Likewise, overview the paper in its first section and 
% summarise in the last.

% \section{Wave scattering model}\label{section:Wave model}
\section{Scatterer representation: splines}\label{Section:Scatterer representation}

\begin{figure}
\centering
\begin{tikzpicture}
\matrix[column sep=10pt, row sep=0pt] {
   \node{\includegraphics[width=0.43\textwidth]{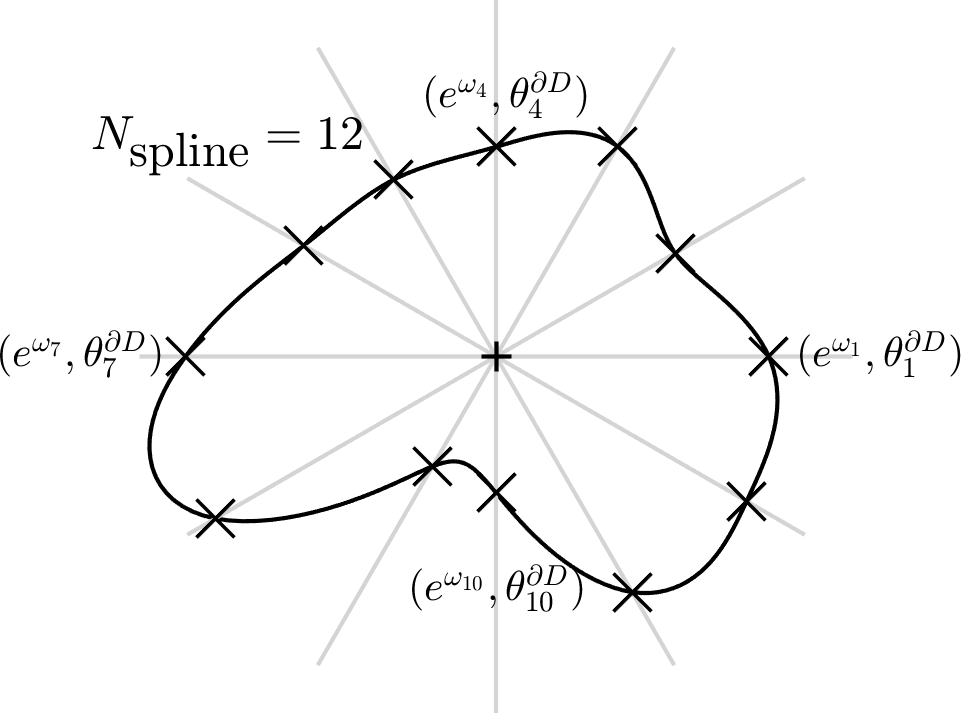}}; &
   \node{\includegraphics[width=0.35\textwidth]{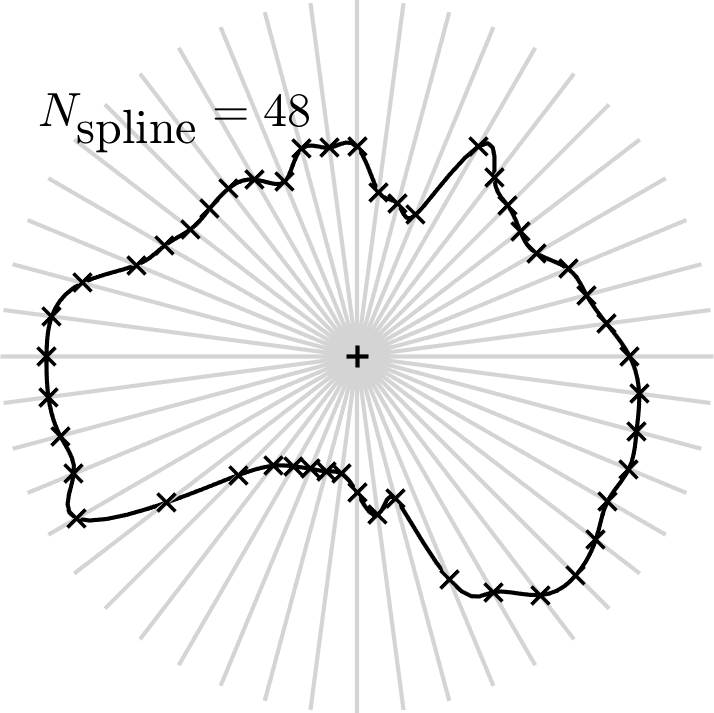}}; \\
};
\end{tikzpicture}
\caption{Star-shaped cubic spline scatterers $D$ approximating the shape in 
Figure \ref{fig:Aus True Cartoon}, for $N_\text{spline}=12$, $48$. The ``$\boldsymbol{+}$'' symbol represents the scatterers center and the ``$\boldsymbol{\times}$'' visualise the knots 
and their associated data. The left plot includes reference $(r,\theta^{\partial D})$ boundary values.}
\label{fig:aus spline approx}
\end{figure}

Let $D(\boldsymbol{\omega}) \subset \mathbb{R}^2$ denote the simply connected region occupied by the scatterer, where $\boldsymbol{\omega}$ denotes a vector characterising the domain $D$. We assume that $D(\boldsymbol{\omega})$ is \textit{star-shaped} with centre at the origin (see Fig \ref{fig:aus spline approx}), and use polar coordinates to parametrise $\partial D (\boldsymbol{\omega})$ through the continuous map $\boldsymbol{\chi}:\left[ 0,2 \pi \right) \rightarrow \partial D (\boldsymbol{\omega})$ given by 
 \begin{equation}
     \boldsymbol{\chi}(\theta) = r(\theta) \boldsymbol{\hat{e}_r}(\theta), \qquad \theta \in [0, 2\pi),
     \label{scatterer polar representation}
 \end{equation}
where $r(\theta)>0$, and $\boldsymbol{\hat{e}_r}(\theta) = \cos \theta \boldsymbol{\hat{e}_x} + \sin \theta \boldsymbol{\hat{e}_y}$ represents the (unit) vector in the radial direction and $\boldsymbol{\hat{e}_x}, \boldsymbol{\hat{e}_y}$ are the constant Cartesian unit vectors.  The closedness of the boundary implies that $\boldsymbol{\chi}(0)=\boldsymbol{\chi}(2 \pi)$. 

We are interested in a class of scatterers in which the \textit{log-radius}
$s(\theta)$,
satisfying
\begin{equation}
    {r}(\theta) =  e^{s(\theta)}, %\quad \theta \in \left[ 0,2\pi \right),
    \label{eq:logradius}
\end{equation}
is a cubic spline
utilising piecewise $C^2$ polynomials and periodic end conditions, so that
\begin{equation}\label{eq:log-radius s(theta)}
s(\theta) = \sum_{j=1}^{N_{\text{spline}}} a_{j}(\boldsymbol{\omega})(\theta - \theta_{j}^{\partial D})^{3} + b_{j}(\boldsymbol{\omega})(\theta - \theta_{j}^{\partial D})^{2} + c_{j}(\boldsymbol{\omega})(\theta - \theta_{j}^{\partial D}) + d_{j}(\boldsymbol{\omega}),
\end{equation}
for $\theta \in [\theta_j^{\partial D},\theta_{j+1}^{\partial D})$. The equispaced knots %satisfy $0= \theta_1^{\partial D} < \theta_2^{\partial D} <... %<\theta_{N_\text{spline}}^{\partial D}$ and 
are $\theta_j^{\partial D} =2\pi (j-1)/N_\text{spline}$ for $j=1,\dots,N_\text{spline}$. 
%{Furthermore, our sensitivity analysis is restricted to the value at the knots %$r(\theta_i^{\partial D}).$ As such, given our assumptions on the interpolating cubics, 
In practice, the values 
$d_i \equiv s(\theta_i^{\partial D})$ for $i=1,\dots,N_\text{spline}$ 
of the log-radius at the knots
are sufficient to uniquely describe the scatterer, and without loss of generality we can write the vector of coefficients $\mathbf{d}(\boldsymbol{\omega}) = \boldsymbol{\omega} \in \mathbb{R}^{N_\text{spline}}$,
and $\mathbf{a}(\boldsymbol{\omega}),\mathbf{b}(\boldsymbol{\omega}),\mathbf{c}(\boldsymbol{\omega})$ in~\eqref{eq:log-radius s(theta)}
are obtained by solving an associated linear system where only the right hand side
depends on $\boldsymbol{\omega}$.
Scatterers obtained
using the construction~\eqref{eq:log-radius s(theta)} for the
shape in Figure~\ref{fig:Aus True Cartoon} with $N_\text{spline}=12,48$
are depicted in Figure \ref{fig:aus spline approx}. 
%In our implementation, we achieve this by using the \texttt{csape} function with the %`periodic' condition in the Matlab (2020a) Curve Fitting Toolbox. 

\section{Scattering model for sensitivity analysis}\label{section: Model}

We consider time-harmonic waves with angular frequency $\psi$ 
propagating in a non-dispersive homogeneous medium with constant wave speed $c$ exterior to 
$D(\boldsymbol{\omega})$. The respective waves are described by the complex valued function $u(\mathbf{x};\boldsymbol{\omega})$, which
satisfies
the 2D Helmholtz equation
\begin{equation}
    \left( \Delta  + k ^2 \right) u(\mathbf{x}; \boldsymbol{\omega}) = 0, \quad \mathbf{x} \in \mathbb{R}^2 \setminus \overline{D}(\boldsymbol{\omega}),
    \label{HH eqn}
\end{equation}
where $k = \psi/c$ is the wavenumber. The waves are forced by an \textit{incident} plane wave
\begin{equation}
    u^{\text{inc}}(\mathbf{x}) = \me^{\mi k \mathbf{x} \cdot \mathbf{\widehat{{d}}}},
\end{equation}
travelling in the direction of the (unit) vector $\mathbf{\widehat{d}}$, 
which induces a \textit{scattered} field $u^{\text{sc}}$. 
%By writing the total field as $ u(\mathbf{x}; \boldsymbol{\omega}) = u^{\text{inc}}(\mathbf{x}) + u^{\text{sc}}(\mathbf{x}; \boldsymbol{\omega})$, linearity implies that $u^{\text{sc}}$ must also satisfy (\ref{HH eqn}). 
The scattered field must additionally satisfy the radiation condition \cite{schot1992eighty}, \begin{equation}
     \frac{\partial u^{\text{sc}}}{\partial \mathbf{x}}(\mathbf{x}; \boldsymbol{\omega}) - \mi k u^{\text{sc}}(\mathbf{x}; \boldsymbol{\omega})  = o\left(\frac{1}{\sqrt{|\mathbf{x}|}}\right), \quad \text{as} \quad |\mathbf{x}| \to \infty, 
    \label{2D radiation condition}
\end{equation}
uniformly with respect to direction $\widehat{\mathbf{x}}  = \mathbf{x}/|\mathbf{x}| = \boldsymbol{\hat{e}_r}(\theta)  \in \partial B$, where $\partial B \subseteq \mathbb{R}^2$ denotes the set of all directions in the unit circle. The radiation condition (\ref{2D radiation condition}) allows us to introduce the far-field pattern ${u}^\infty (\widehat{\mathbf{x}}; \boldsymbol{\omega})$ induced as a result of the interaction of an incident plane wave with the particle $D(\boldsymbol{\omega})$ via
\begin{equation}
{u}^{\text{sc}}(\mathbf{x}; \boldsymbol{\omega}) = \frac{\me^{\mi k |\mathbf{x}|}}{\sqrt{|\mathbf{x}|}} \left({u}^\infty (\widehat{\mathbf{x}}; \boldsymbol{\omega}) + O \left( \frac{1}{|\mathbf{x}|} \right) \right), \quad \text{as} \quad |\mathbf{x}| \to \infty,
\label{definiton u infty}
\end{equation}
see e.g. (eqn $3.109$ in \cite{colton2019inverse}). 

This work is motivated by the inverse problem in which far-field data is used to reconstruct the shape of the scatterer.
For a fixed incident plane wave,
it is convenient to 
introduce the {far-field} operator $\mathcal{F}: C([0,2\pi)) \rightarrow C(\partial B)$ by 
\begin{equation}
    \mathcal{F}\left[r(\theta, \boldsymbol{\omega})\right] =  u^\infty(\widehat{\mathbf{x}};\boldsymbol{\omega}), \label{F(omega)}
\end{equation} 
which maps the boundary $\partial D$ (defined using $r$) of the sound--soft scatterer onto the far-field pattern of the scattered wave. For a sound-soft
scatterer with Dirichlet boundary conditions on the boundary $\partial D(\boldsymbol{\omega})$, we have
\begin{equation}
{u}(\mathbf{x}; \boldsymbol{\omega}) = f, \quad \text{for} \quad \mathbf{x} \in \partial D(\boldsymbol{\omega}),
\label{Dirichlet BC}    
\end{equation}
with $f=0$.

An efficient way to evaluate
$\mathcal{F}$ when $f=0$ is to use the BIE reformulation
\begin{equation}\label{burton miller}
    \left(\mathcal{I} + 2 \mathcal{K}' - \mi k\mathcal{S} \right) \left[ \frac{\partial u}{\partial \mathbf{n}} \right] = 2\frac{\partial u^\text{inc}}{\partial \mathbf{n}}  -\mi k u^\text{inc} , \quad \mathbf{x} \in \partial D(\boldsymbol{\omega}),
\end{equation}
where $\mathcal{I}$ is the identity operator, $\mathcal{S}$ is the single-layer potential, and $\mathcal{K}'$ is the
double-layer transpose potential~\cite{colton2019inverse}.
The unknown in~\eqref{burton miller} is the physical-quantity
$\partial u/\partial \mathbf{n}$,  where $\mathbf{n}(\mathbf{x_0})$ denotes the outward unit normal to $\partial D(\boldsymbol{\omega})$ at $\mathbf{x_0}$).
Once $\partial u/\partial \mathbf{n}$ has been computed by solving~\eqref{burton  miller},
we compute the scattered field away from the boundary as $u^\text{sc}= -\mathcal{S} [\partial u/\partial \mathbf{n}]$ for $\mathbf{x} \in \mathbb{R}^2 \setminus \overline{D}(\boldsymbol{\omega})$.
The far-field can then be computed similarly.

We compute high-order approximate solutions to the wave-scattering problem introduced above using a spectral {Nystr\"om} method, implemented in open source software \cite{ganesh2017algorithm}, in which equation (\ref{burton miller}) is discretized using equispaced points on the boundary $\partial D(\boldsymbol{\omega})$, leading to a linear system 
% $\mathbf{A}(\boldsymbol{\omega}) \boldsymbol{\phi}(\boldsymbol{\omega}) = \mathbf{b}$,
\begin{equation}\label{eqn: nystrom matrix system}
    \mathbf{A}(\boldsymbol{\omega}) \boldsymbol{\phi}(\boldsymbol{\omega}) = \mathbf{b},
\end{equation}
where the matrix $\mathbf{A}(\boldsymbol{\omega})$ represents the discretization of the combined field integral operators on the left hand side of (\ref{burton miller}), $\mathbf{b}$ contains samples of the right hand side of (\ref{burton miller}) evaluated at the discretization points, and $\boldsymbol{\phi}$ approximates $\partial u/\partial \mathbf{n}$ at these points. %The $C^2$ continuity requirement on the spline representation of $\partial D(\boldsymbol{\omega})$ discussed in Section \ref{subSection:Scatterer representation} ensures suitable regularity for high-order quadrature.
We emphasise that, because the integration in~\eqref{burton miller} is over 
the boundary $\partial D$, the matrix $\mathbf{A}(\boldsymbol{\omega})$ must be recomputed for each new scatterer configuration. 
%After solving for $\boldsymbol{\phi}$ in (\ref{eqn: nystrom matrix system}), the far-field pattern $u^\infty$ can be evaluated efficiently using the same quadrature scheme.

%\section{Sensitivity analysis}\label{sec: shape der}
Next we examine the sensitivity of the far field to changes in the scatterer's data, by
considering the partial derivatives of the 
mapping $\mathcal{F}:r(\theta;\boldsymbol{\omega}) \mapsto u^\infty$.
Using the chain rule gives
\begin{equation}\label{dF/d omega_j}
    J_i = \frac{\partial  \mathcal{F}}{\partial \omega_i} = \frac{\mathrm{d}  \mathcal{F}}{\mathrm{d} r} \left[ \frac{\partial r}{\partial \omega_i}\right], \qquad \frac{\partial r}{\partial \omega_i} = r(\theta) \frac{\partial s(\theta; \boldsymbol{\omega})}{\partial \omega_i}.
\end{equation}
Here
$\mathrm{d} \mathcal{F}[\cdot]/\mathrm{d} r $ is the \textit{Fr\'echet derivative} 
of $\mathcal{F}$. The Fr\'echet derivative is given by the far-field of an auxiliary wave field $v$ satisfying (\ref{HH eqn}), (\ref{2D radiation condition}) and the inhomogeneous Dirichlet condition (see {\cite[Theorem 5.15]{colton2019inverse}})
\begin{equation}\label{thm: frechet der}
%         \frac{\mathrm{d} \mathcal{F}}{\mathrm{d}r} \left[q\right] = v_\infty, \qquad  
v(\mathbf{x};\boldsymbol{\omega}) = -\mathbf{n} \cdot (q(\theta) \boldsymbol{\hat{e}_r}(\theta) \circ \boldsymbol{\chi}^{-1}) \frac{\partial u}{\partial \mathbf{n}} \quad \text{for}\quad \mathbf{x} \in \partial D(\boldsymbol{\omega}),
\end{equation}
with $q \in \partial B$ a scalar function representing the direction. To compute $v$
we solve the combined-field BIE \cite[(3.29)]{colton2019inverse} for $\phi^*$ where the inhomogeneity is coupled with (\ref{burton miller}),
\begin{equation}\label{eqn: inhomog. dirichlet}
    (\mathcal{I} + 2\mathcal{K} + 2\mi k \mathcal{S})[\phi^*] = - 2 \mathbf{n} \frac{\partial u}{\partial \mathbf{n}} \cdot (q(\theta) \boldsymbol{\hat{e}_r}(\theta) \circ \boldsymbol{\chi}^{-1}), \quad \mathbf{x} \in \partial D(\boldsymbol{\omega}).
\end{equation}
 Once the surface potential $\phi^*$ has been obtained by solving~\eqref{eqn: inhomog. dirichlet}, we compute the induced field using
$v(\mathbf{x}) = (\mathcal{K} + \mi k\mathcal{S})[\phi^*]$ for $\mathbf{x} \in \mathbb{R}^2 \setminus \overline{D}(\boldsymbol{\omega})$.
The far field is computed similarly.

%\begin{remark}
\textit{Remark} 1.
    We use the BIE~\eqref{burton miller} to solve the forward scattering problem
    because it yields the surface potential
$\partial u/\partial \mathbf{n}$, which is required in the right hand side 
of~\eqref{eqn: inhomog. dirichlet} to compute the Fr\'echet derivative. 
The BIE~\eqref{burton miller} is not appropriate for solving the Fr\'echet-derivative
PDE subject to
the Dirichlet data in~\eqref{thm: frechet der}, because the right hand side
of~\eqref{burton miller}
would require a further normal derivative of $u$; in this case the simpler
indirect BIE~\eqref{eqn: inhomog. dirichlet} is preferred.
%    \end{remark}

To better understand the sensitivity of the derivative to perturbations in the shape,
we compute the 
singular value decomposition
\begin{equation}\label{SVD decomposed Jacobian}
    |\mathbf{J}| = \mathbf{U}\mathbf{\Sigma}\mathbf{V}^\top,
\end{equation}
where $J_{ij} = J_i(\theta_j)$ are the entries in the Jacobian matrix,
$\mathbf{\Sigma}$
is the matrix of singular values $\sigma_i \geq 0$, $\mathbf{U} = [\mathbf{u_1}|\mathbf{u_2}|...|\mathbf{u_{N_\text{spline}}}]$, with $\mathbf{u_i}$ the unit left singular vectors and $\mathbf{V} = [\mathbf{v_1}|\mathbf{v_2}|...|\mathbf{v_{N_\text{obs}}}]$ with $\mathbf{v_i}$ the (unit) right singular vectors of $|\mathbf{J}|$. 
For multiple incident waves 
we define $\mathcal{F}$ to be the concatenation of the mappings for each incident wave.
Then an analogue of
(\ref{SVD decomposed Jacobian}) applies and $|\mathbf{J}| \in \mathbb{R}^{N_\text{spline} \times (N_{\text{inc}}N_{\text{obs}})}$ is the reshaped sensitivity matrix including all incident directions. This is useful since, as we will see below, it deemphasises the role of the incident wave direction in identifying the most significant geometric features contributing to the far-field.

\begin{figure}
\centering
\begin{tikzpicture}
\matrix[column sep=1pt, row sep=0pt] {
&\node{\includegraphics[width=0.3\textwidth]{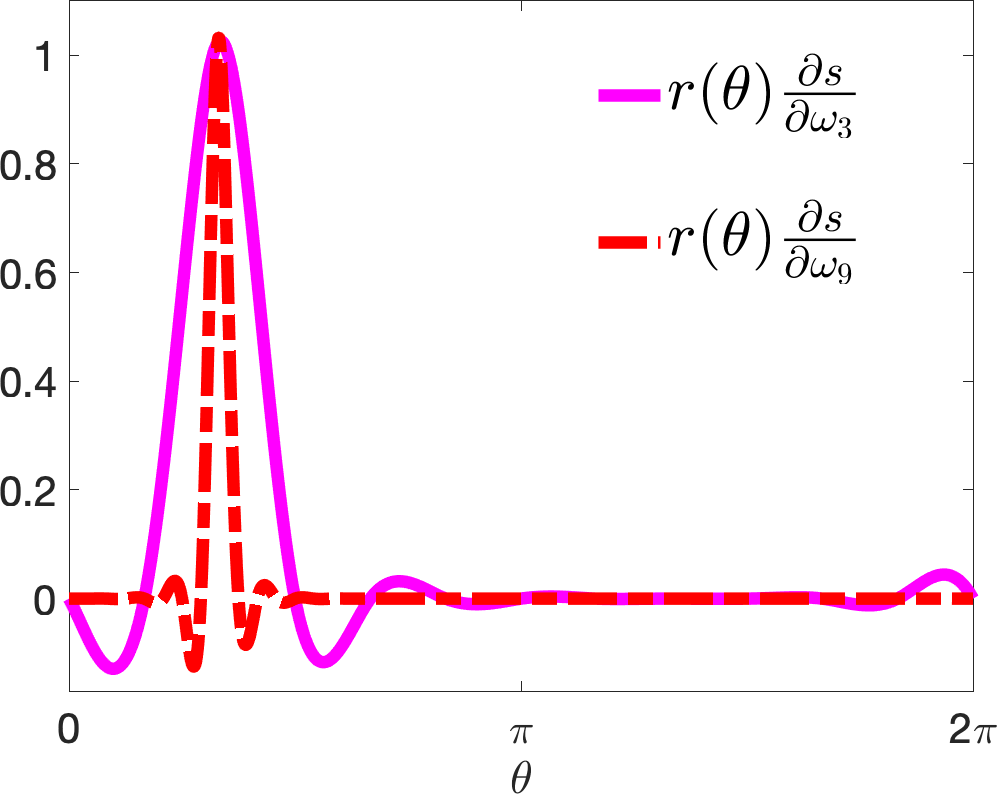}}; &
 &\node{\includegraphics[width=0.3\textwidth]{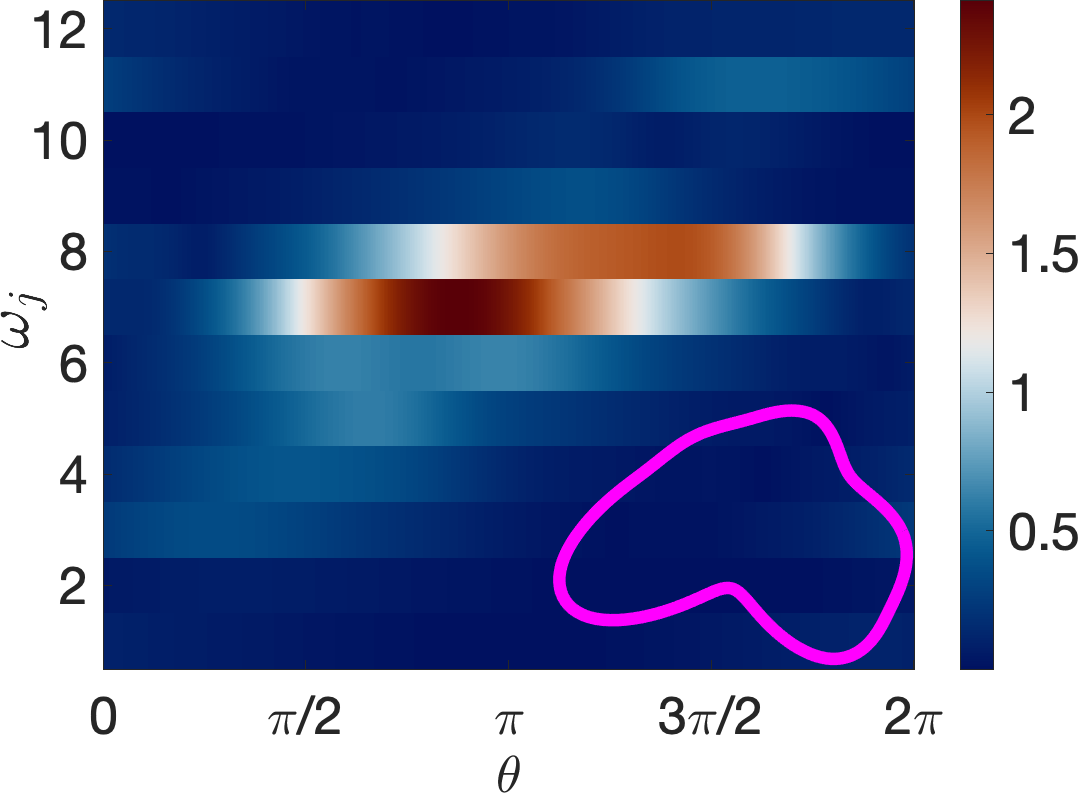}}; &
   \node{\includegraphics[width=0.3\textwidth]{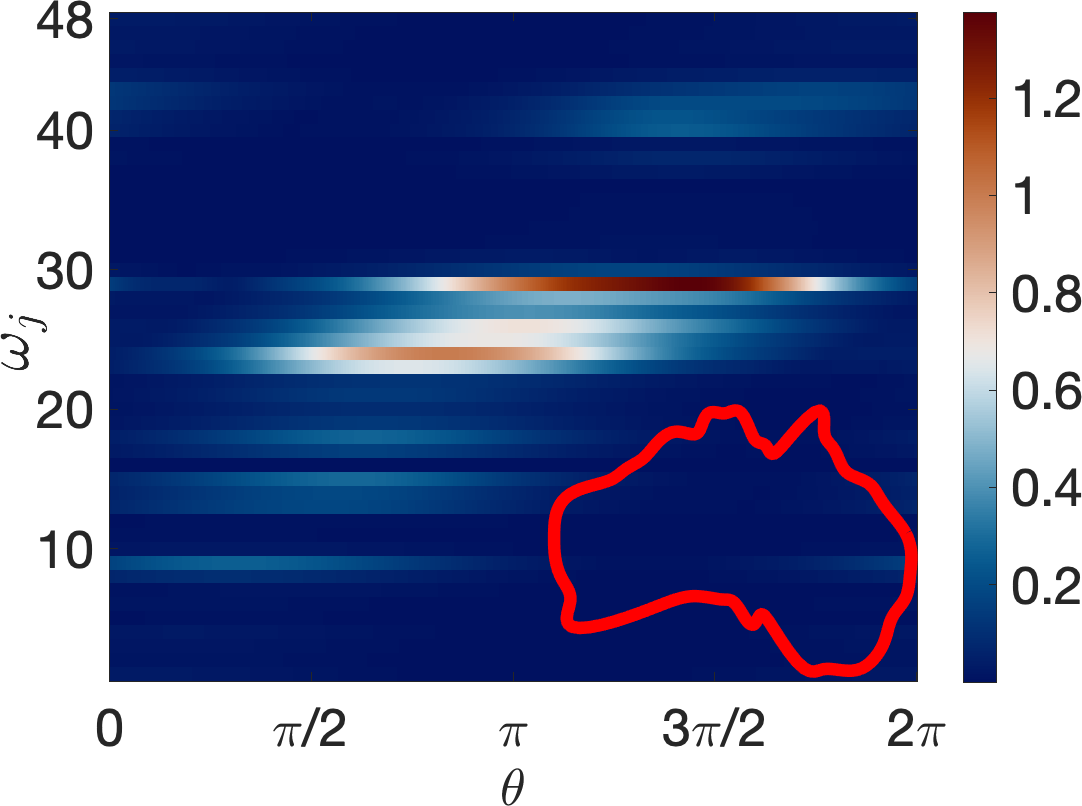}};\\
   };
\end{tikzpicture}
\caption{Left: Input function to the Fr\'echet operator (\ref{dF/d omega_j}) for our 2 scatterers at $\theta_i^{\partial D}=\pi/3$. Jacobian matrix $|\mathbf{J}|$ for $\mathbf{\widehat{{d}}}=\boldsymbol{\hat{e}_x}$ for $ka=2\pi$, $N_\text{spline}=12$ (center) and $N_\text{spline}=48$ (right).}
\label{fig:1 inc wave matrix bands}
\end{figure}

\section{Numerical results}\label{section:results}
We present results for the scatterers in Figure \ref{fig:aus spline approx} with $N_\text{spline}=12,48$. 
From~\eqref{eq:log-radius s(theta)}, it can be seen that the 
the input function 
$q = \partial r/\partial \omega_i$
to the Fr\'echet operator (\ref{dF/d omega_j}) 
is itself a cubic spline whose values at the knots are
$\partial s(\theta_i^{\partial D};\boldsymbol{\omega})/\partial \omega_i = \mathbf{e}_i$.
It follows that the associated perturbation of the shape is localised around $\omega_i$,
as is illustrated in Figure \ref{fig:1 inc wave matrix bands} (left).

We let $a$ be the diameter of the smallest ball circumscribing the scatterer
and present our results in terms of the non-dimensional frequency $ka$.
Figure \ref{fig:1 inc wave matrix bands} (center, right) visualise the Jacobian matrix 
for a single incident wave
with direction $\mathbf{\widehat{{d}}}=\boldsymbol{\hat{e}_x}$ and $ka=2 \pi$. Here we use $N_\text{obs}=1000$ and
${N_\text{Nys}}=27,240$ for $N_\text{spline}=12,48$ respectively, with the higher
order quadrature required in the latter case to resolve the more complex geometry. 
These Jacobian plots contain key physical information, but the interplay between incident direction, observation angle, and geometric features makes interpretation challenging. We focus on these effects below making use of (\ref{SVD decomposed Jacobian}).

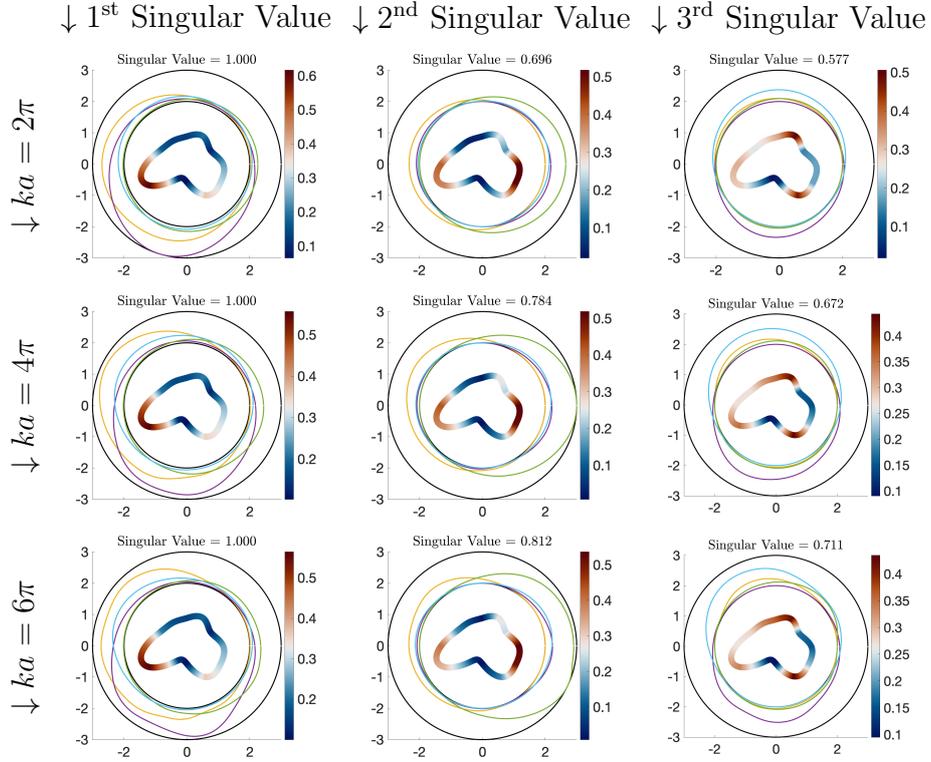
\begin{figure}
\centering
\begin{tikzpicture}
\matrix[column sep=0pt, row sep=0pt] {
    &
   \node  {$\downarrow 1^\text{st}$ \text{Singular Value}}; & \node  {$\downarrow 2^\text{nd}$ \text{Singular Value}};
   & \node  {$\downarrow 3^\text{rd}$ \text{Singular Value}};\\
   \node{\rotatebox{90}{$\downarrow ka=2 \pi$}}; &\node{\includegraphics[width=0.23\textwidth]{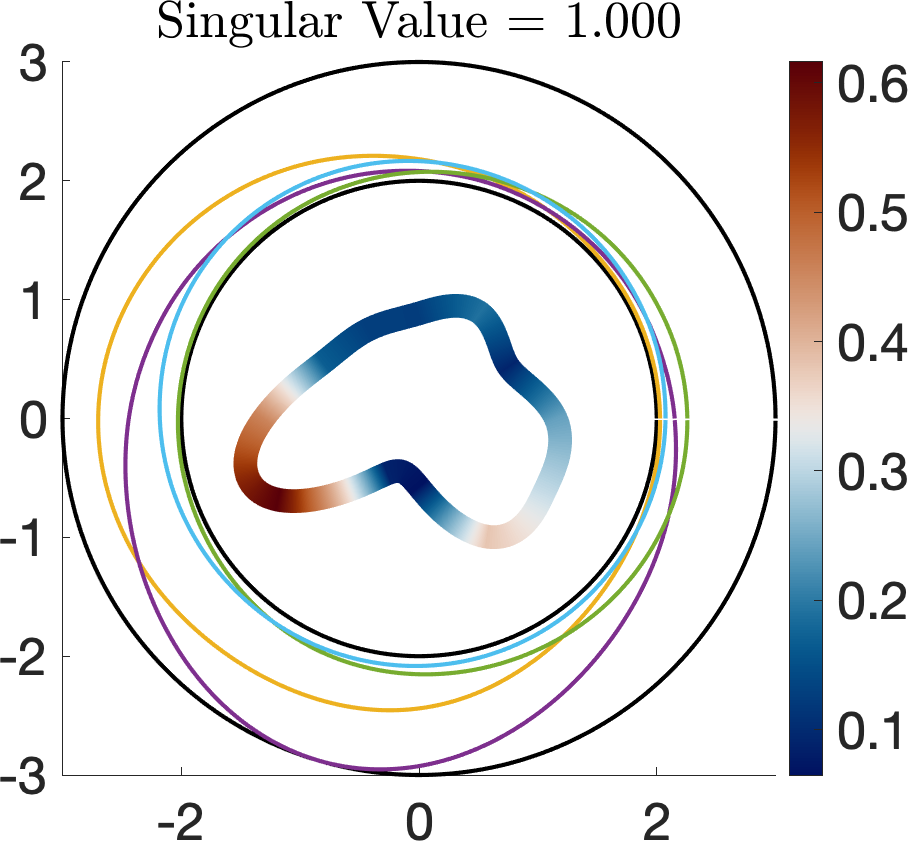}}; &
   \node{\includegraphics[width=0.23\textwidth]{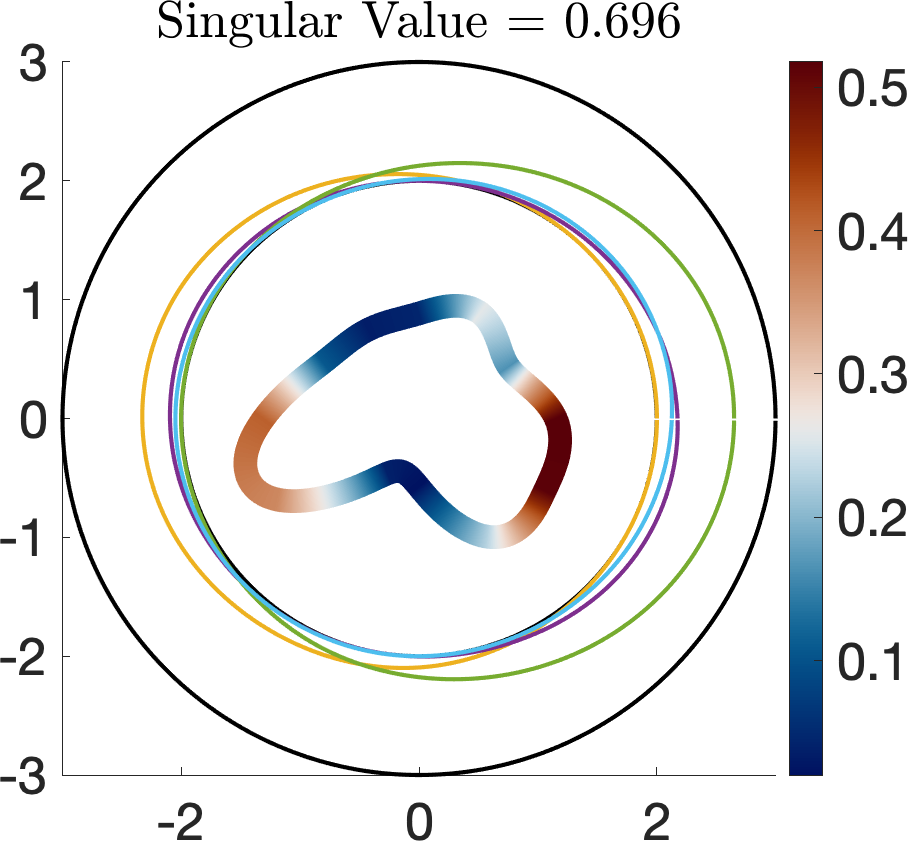}}; &
   \node{\includegraphics[width=0.23\textwidth]{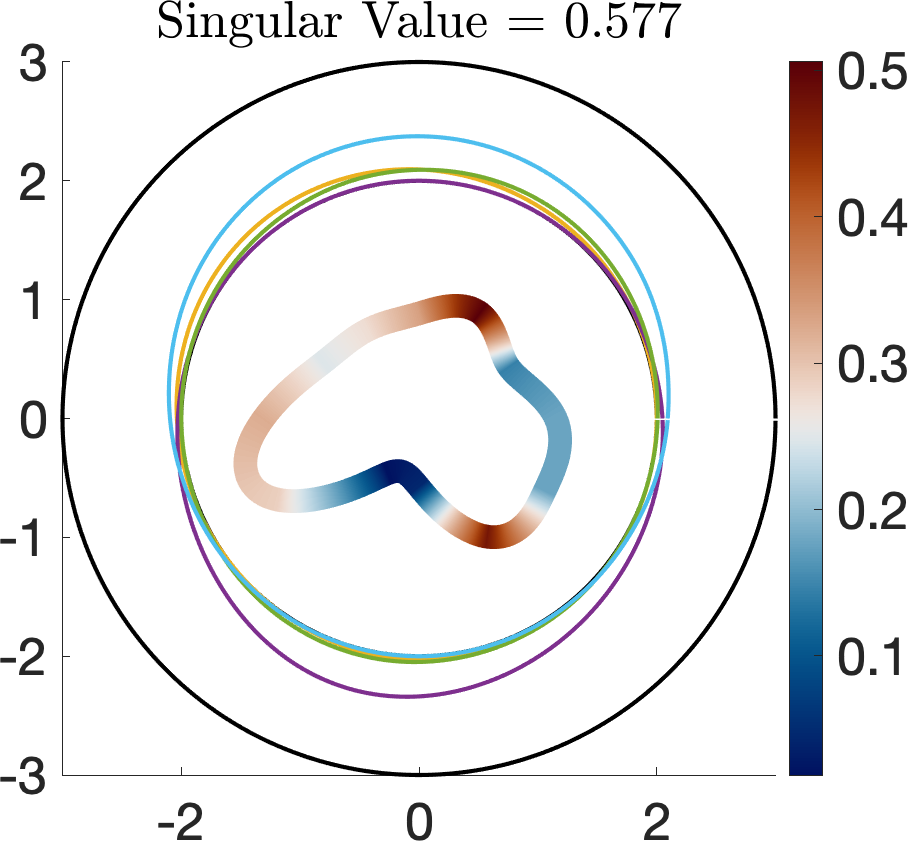}}; \\
    \node{\rotatebox{90}{$\downarrow ka=4 \pi$}};  &\node{\includegraphics[width=0.23\textwidth]{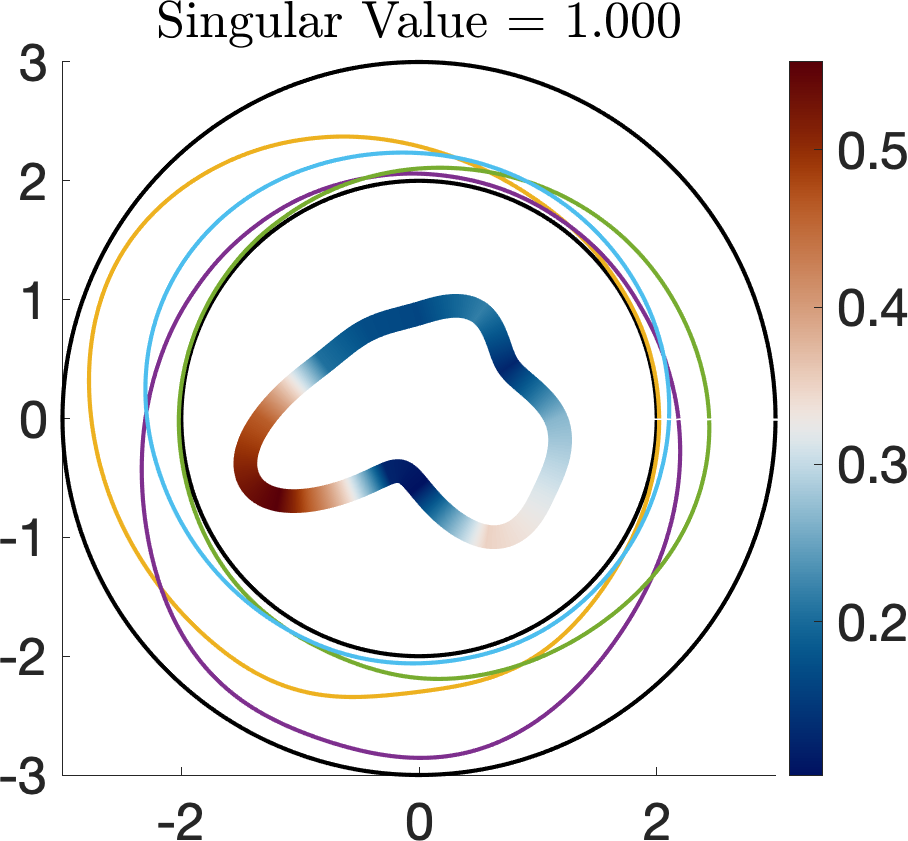}}; &
   \node{\includegraphics[width=0.23\textwidth]{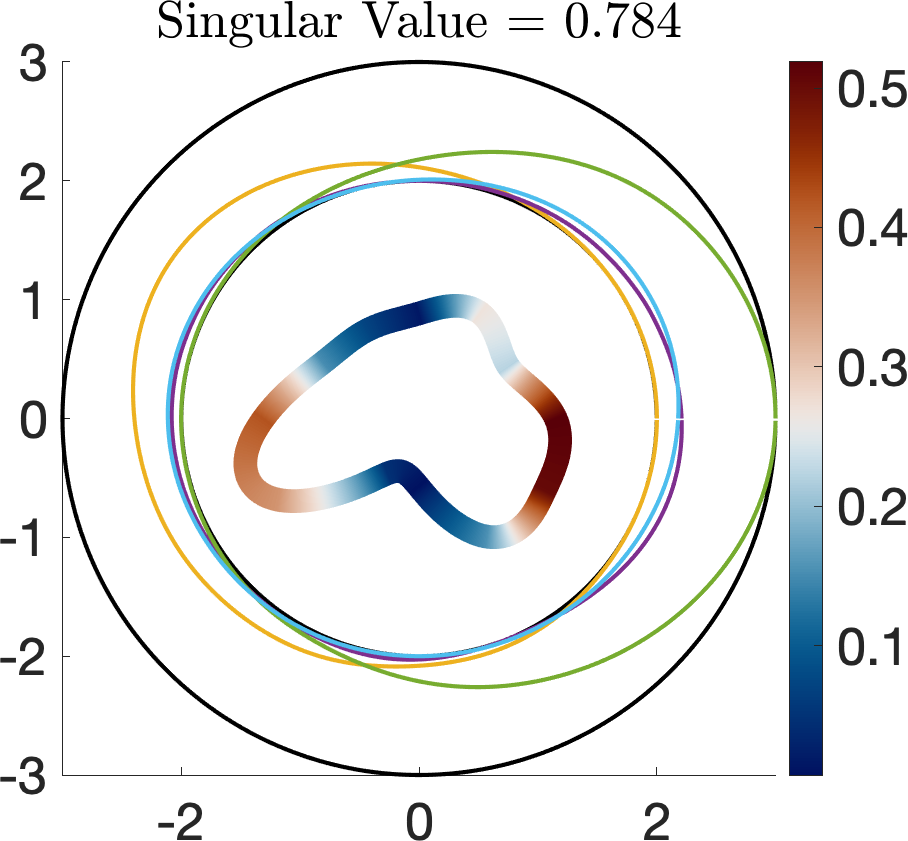}}; &
   \node{\includegraphics[width=0.23\textwidth]{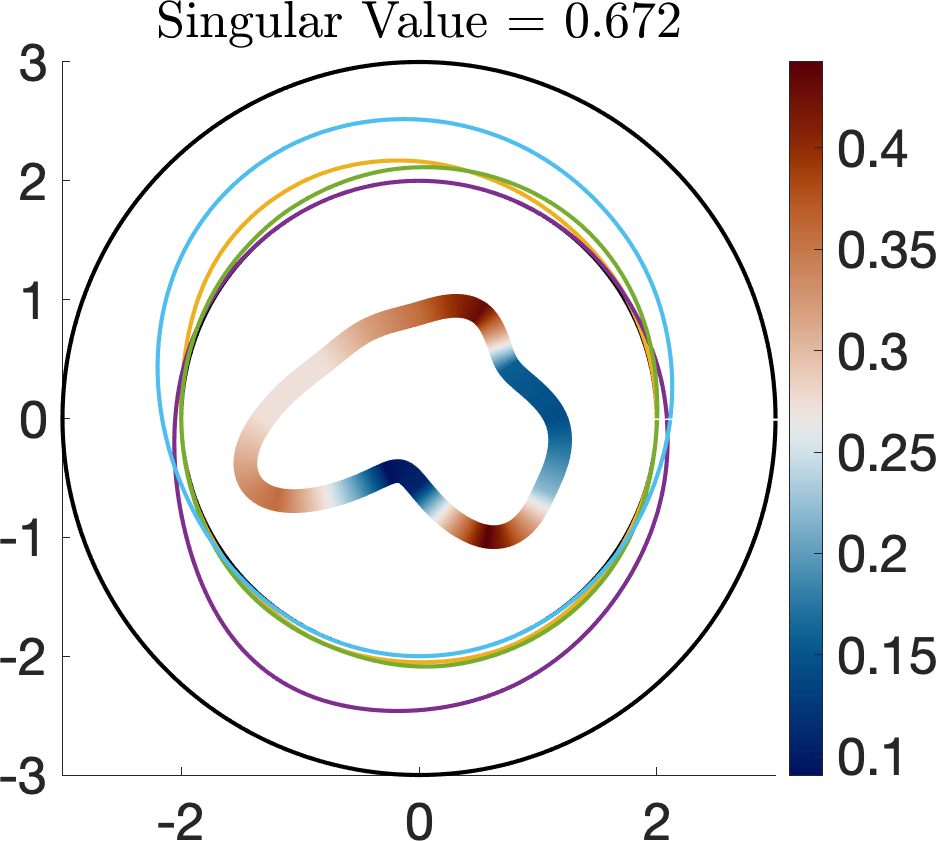}}; \\
    \node{\rotatebox{90}{$\downarrow ka=6 \pi$}}; &\node{\includegraphics[width=0.23\textwidth]{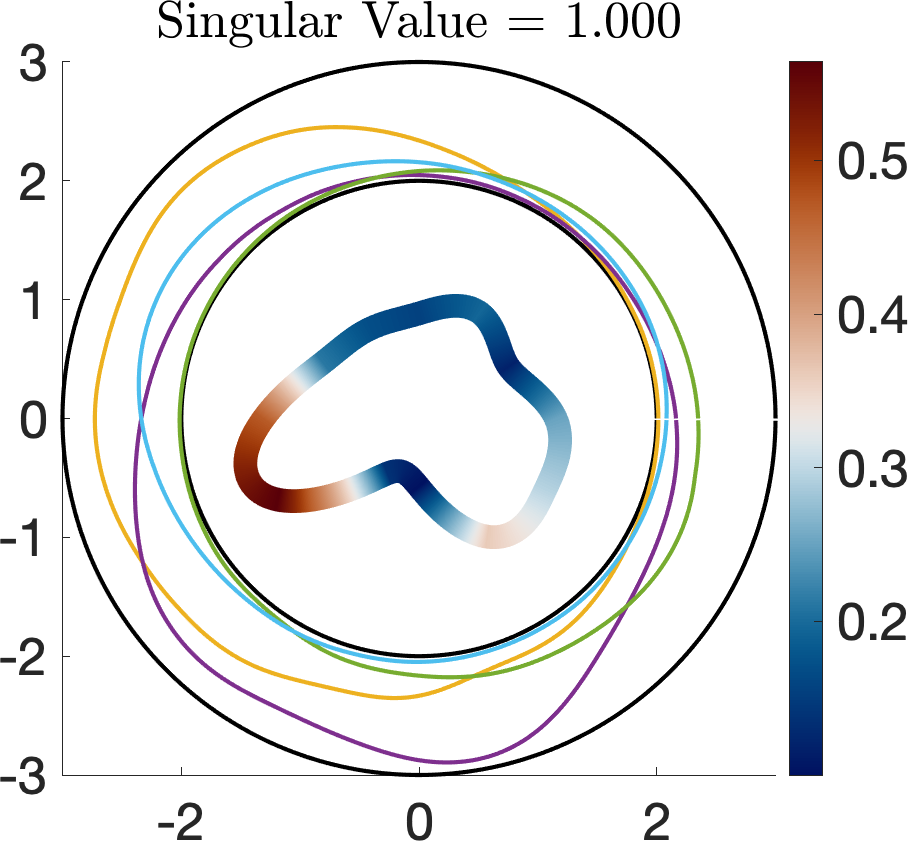}}; &
   \node{\includegraphics[width=0.23\textwidth]{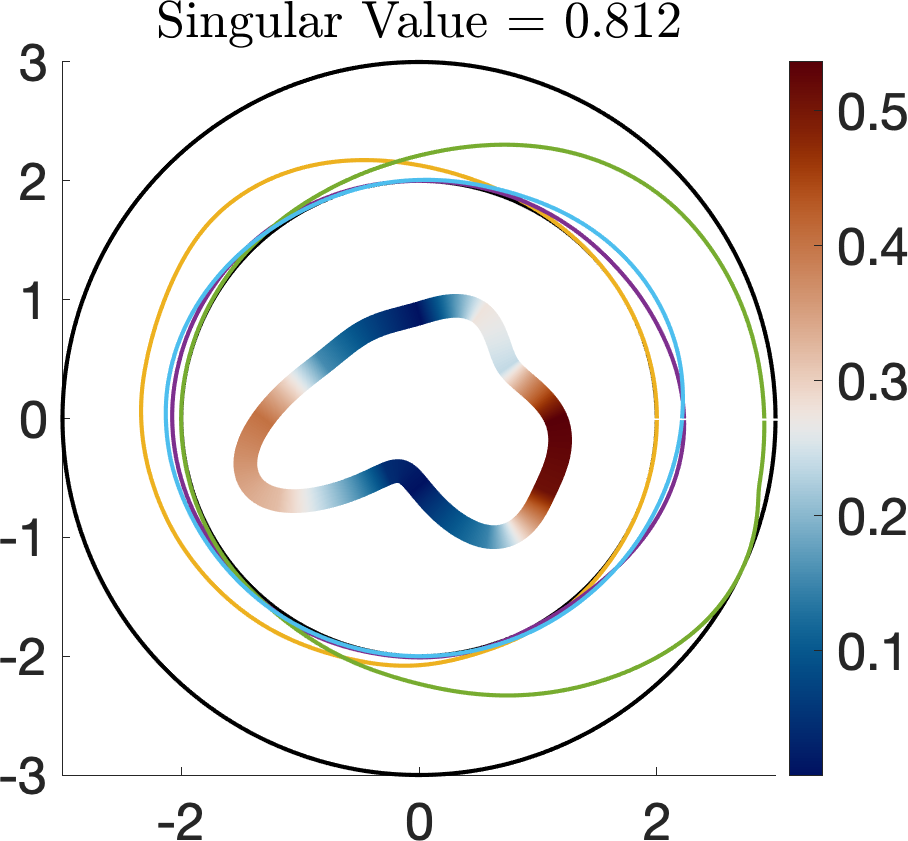}}; &
   \node{\includegraphics[width=0.23\textwidth]{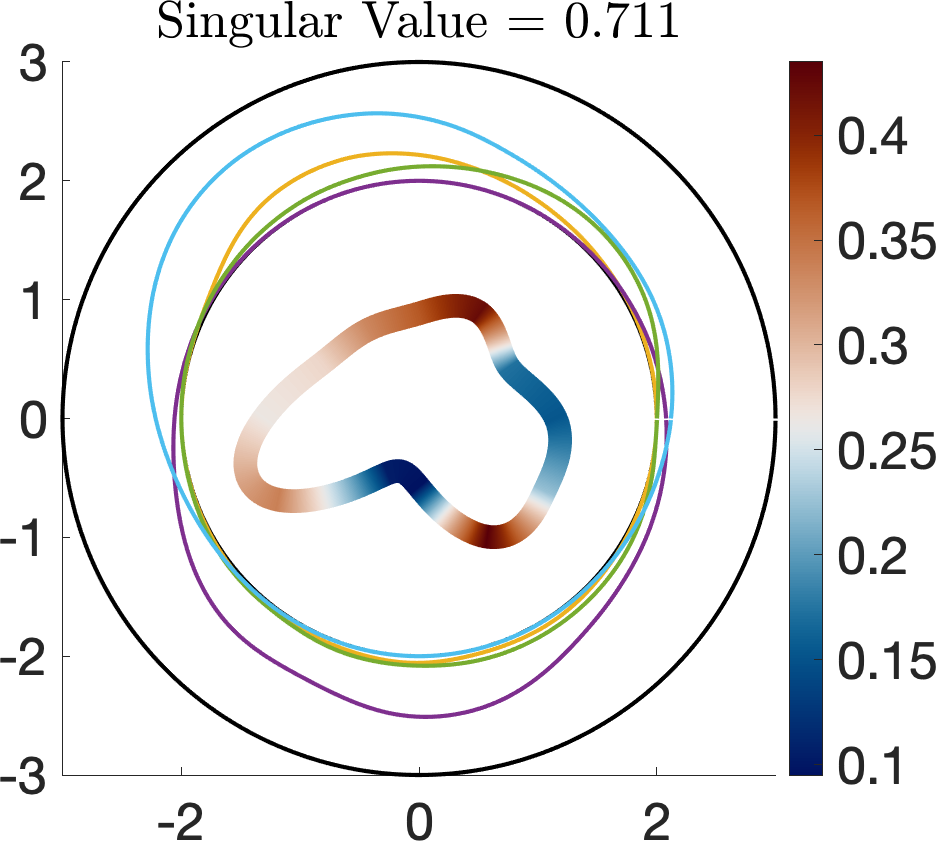}}; \\
   };
\end{tikzpicture}
\caption{Visualisation of the first 3 singular values and associated right and left singular vectors for $N_\text{spline}=12$ and $|\mathbf{J}|$ comprising of $4$ incident waves with $\mathbf{\widehat{{d}}}= \pm \boldsymbol{\hat{e}_x}, \pm \boldsymbol{\hat{e}_y}$.}
\label{fig:4 inc wave SVD N=12}
\end{figure}

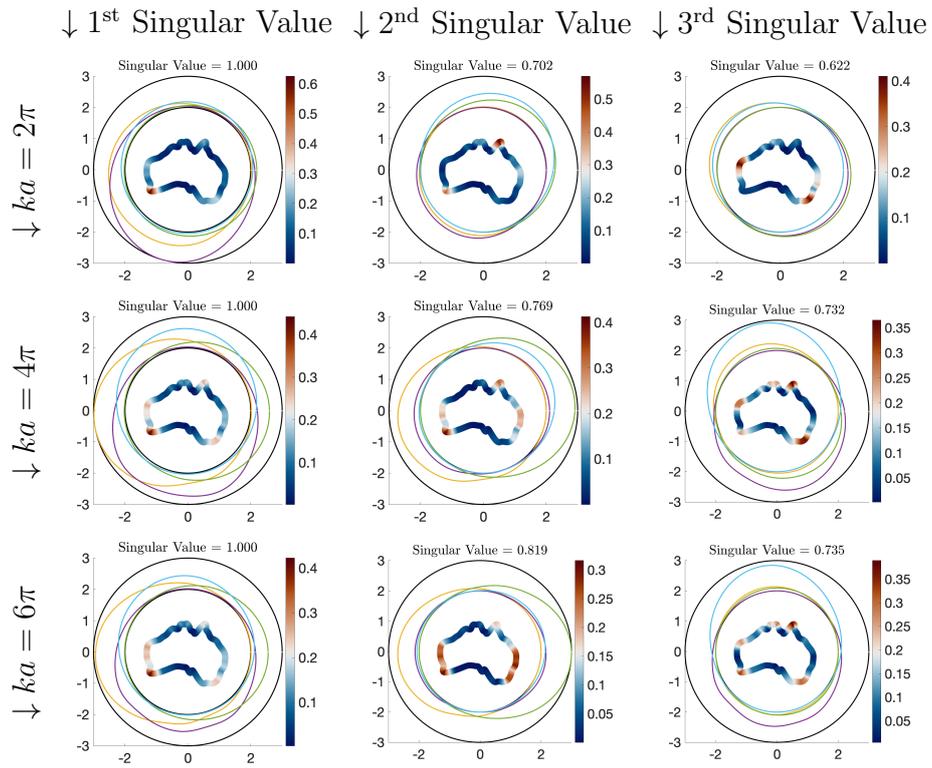
\begin{figure}
\centering
\begin{tikzpicture}
\matrix[column sep=0pt, row sep=0pt] {
    &
   \node  {$\downarrow 1^\text{st}$ \text{Singular Value}}; & \node  {$\downarrow 2^\text{nd}$ \text{Singular Value}};
   & \node  {$\downarrow 3^\text{rd}$ \text{Singular Value}};\\
   \node{\rotatebox{90}{$\downarrow ka=2 \pi$}}; &\node{\includegraphics[width=0.23\textwidth]{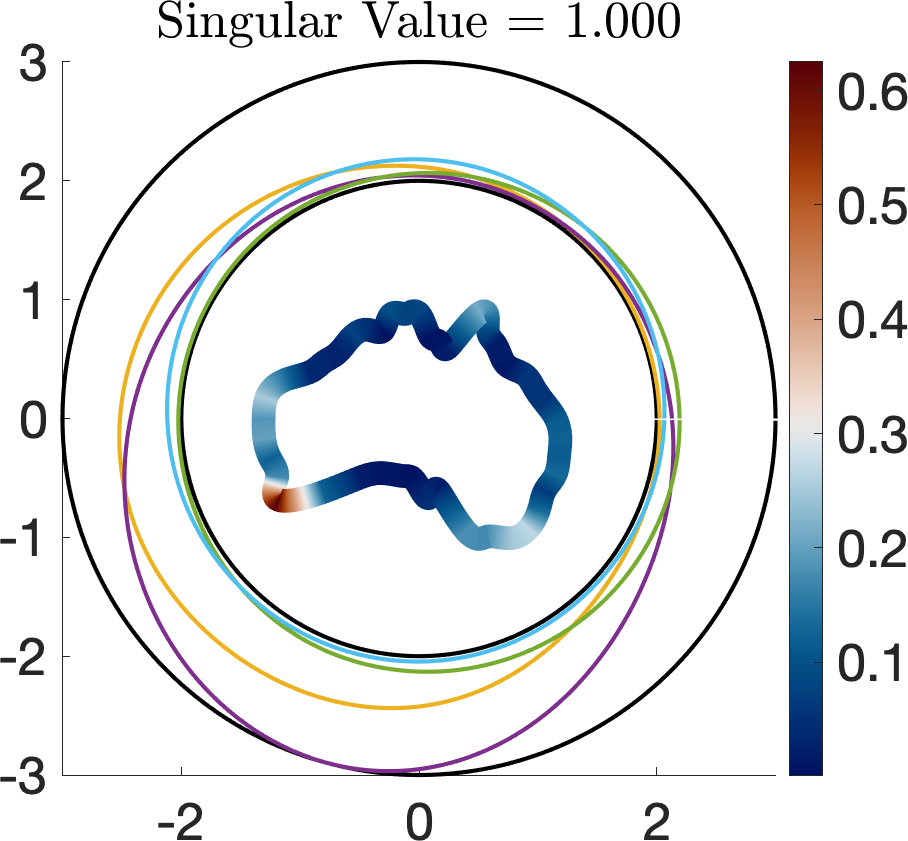}}; &
   \node{\includegraphics[width=0.23\textwidth]{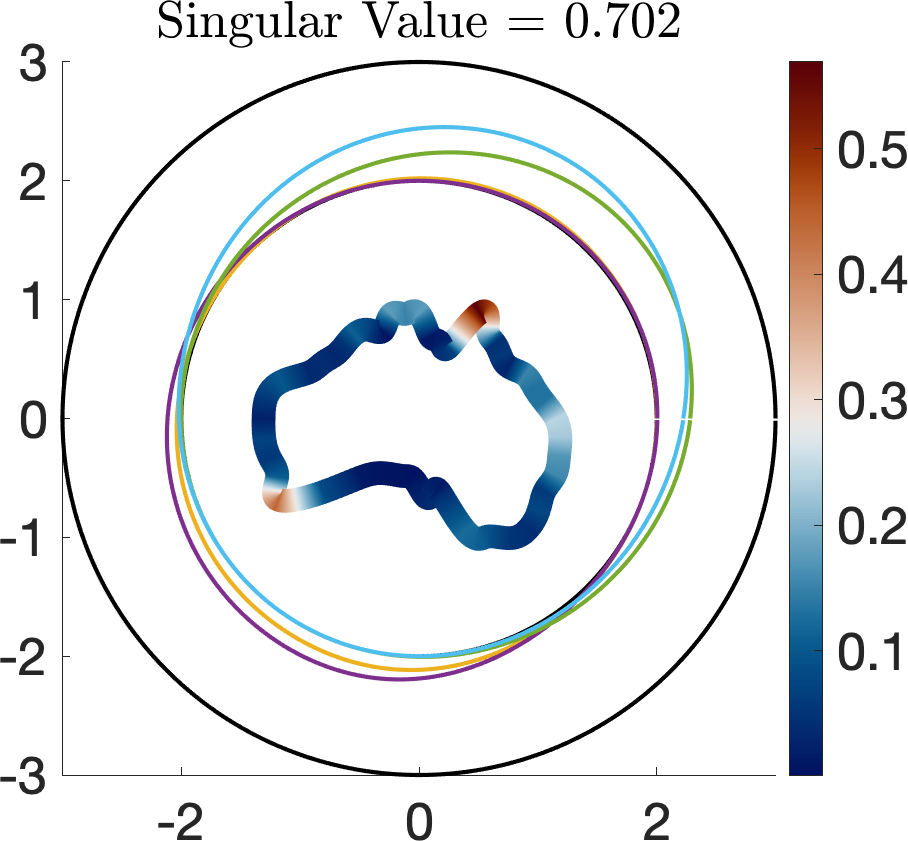}}; &
   \node{\includegraphics[width=0.23\textwidth]{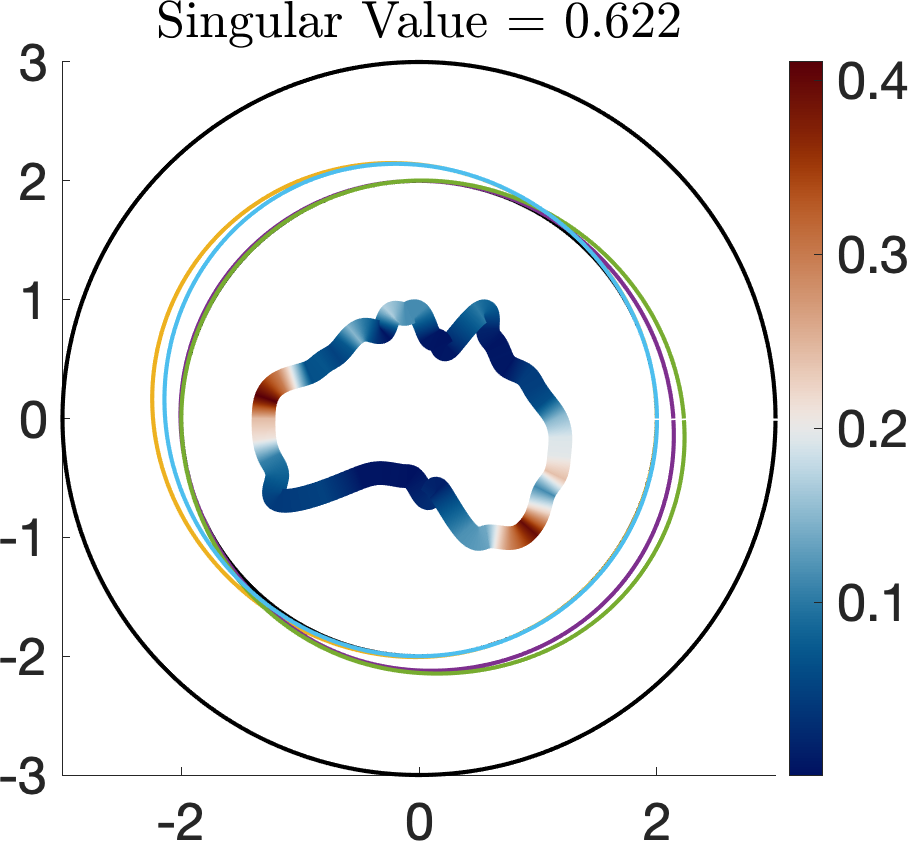}}; \\
    \node{\rotatebox{90}{$\downarrow ka=4 \pi$}};  &\node{\includegraphics[width=0.23\textwidth]{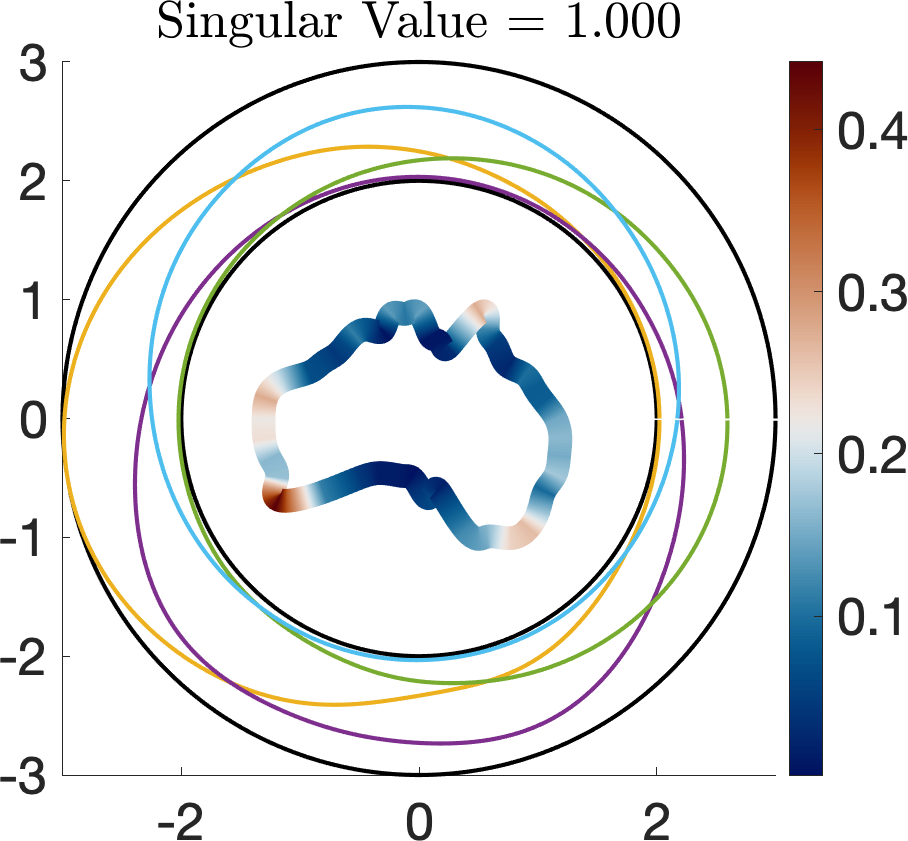}}; &
   \node{\includegraphics[width=0.23\textwidth]{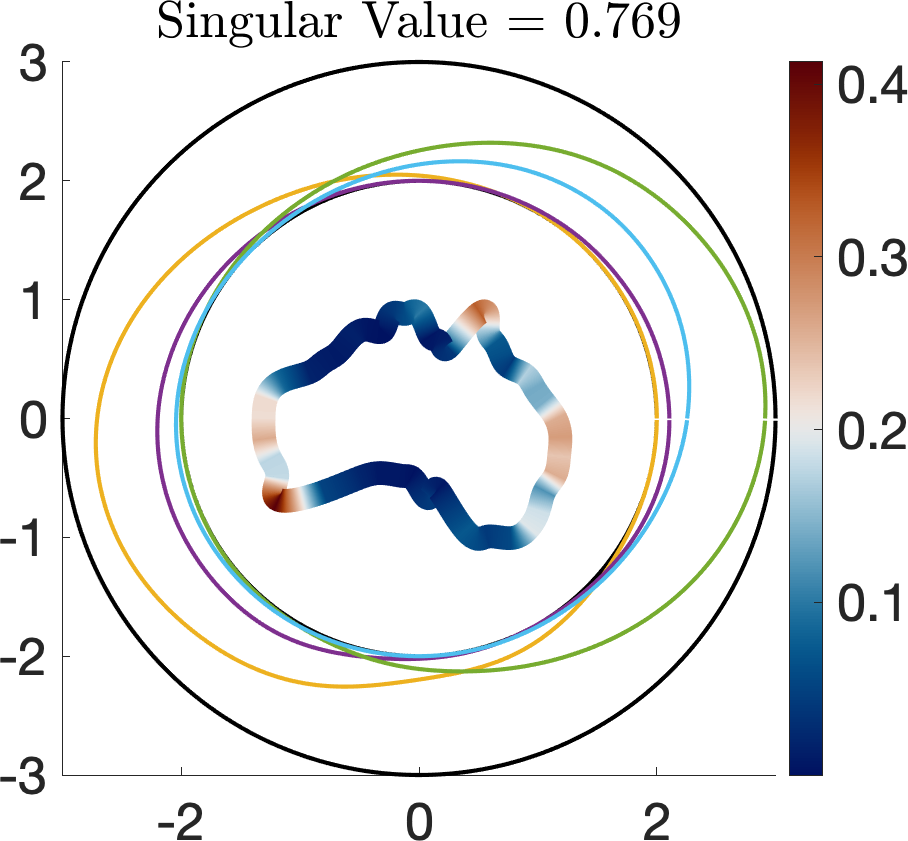}}; &
   \node{\includegraphics[width=0.23\textwidth]{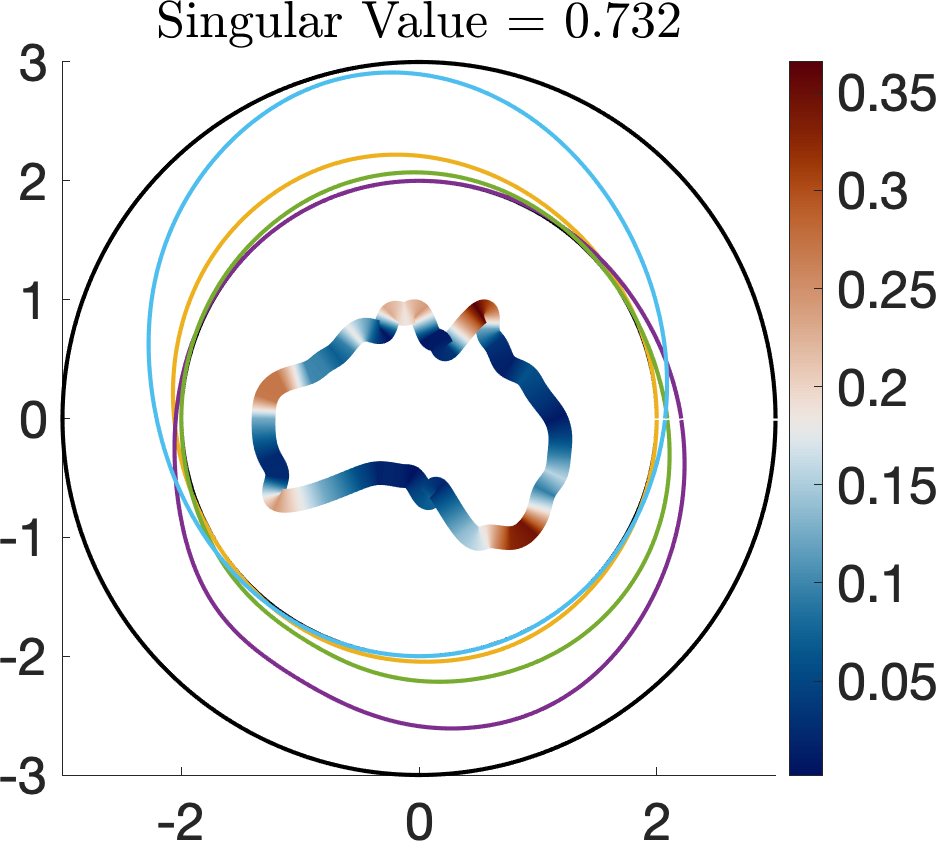}}; \\
    \node{\rotatebox{90}{$\downarrow ka=6 \pi$}}; &\node{\includegraphics[width=0.23\textwidth]{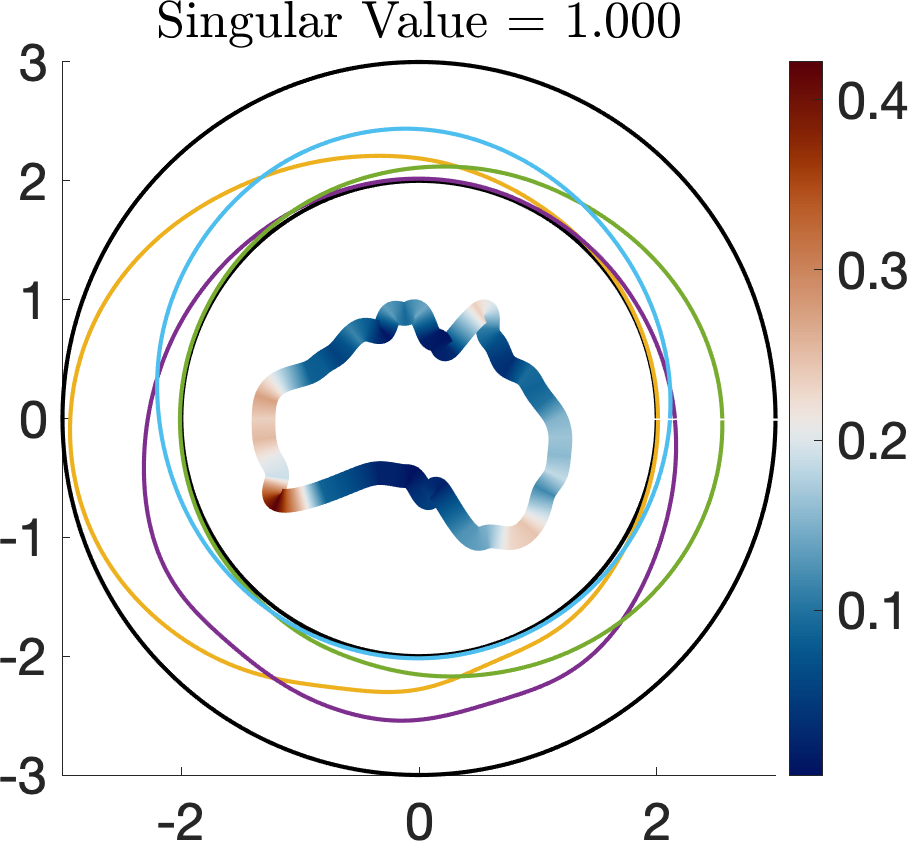}}; &
   \node{\includegraphics[width=0.23\textwidth]{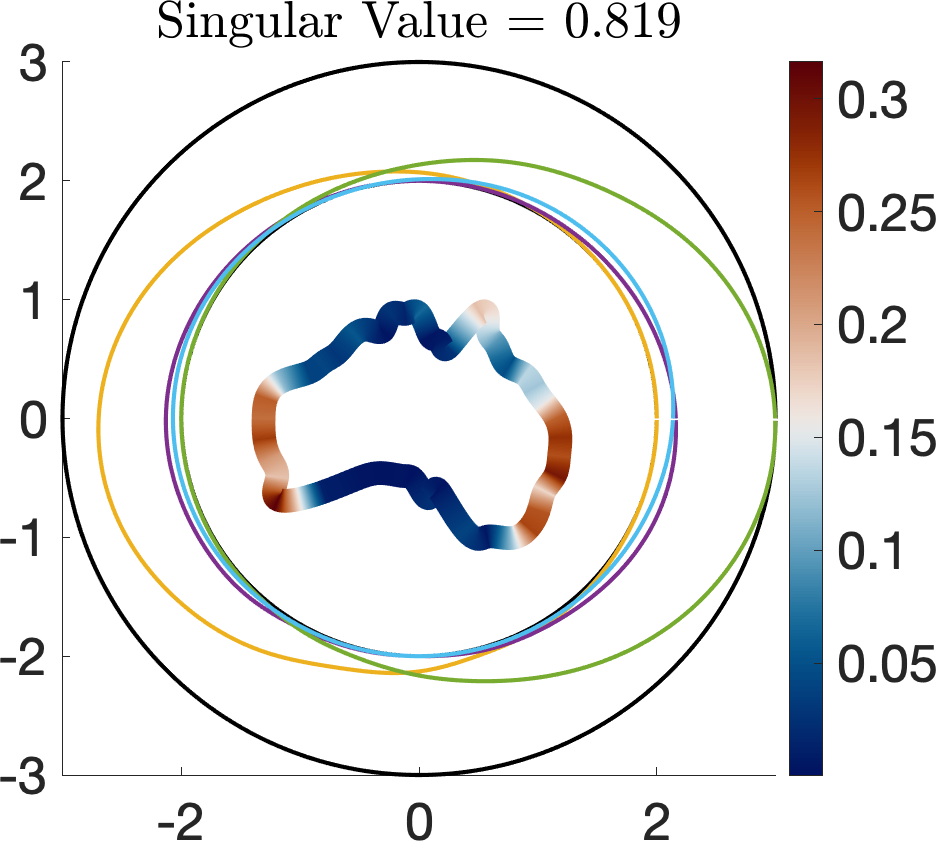}}; &
   \node{\includegraphics[width=0.23\textwidth]{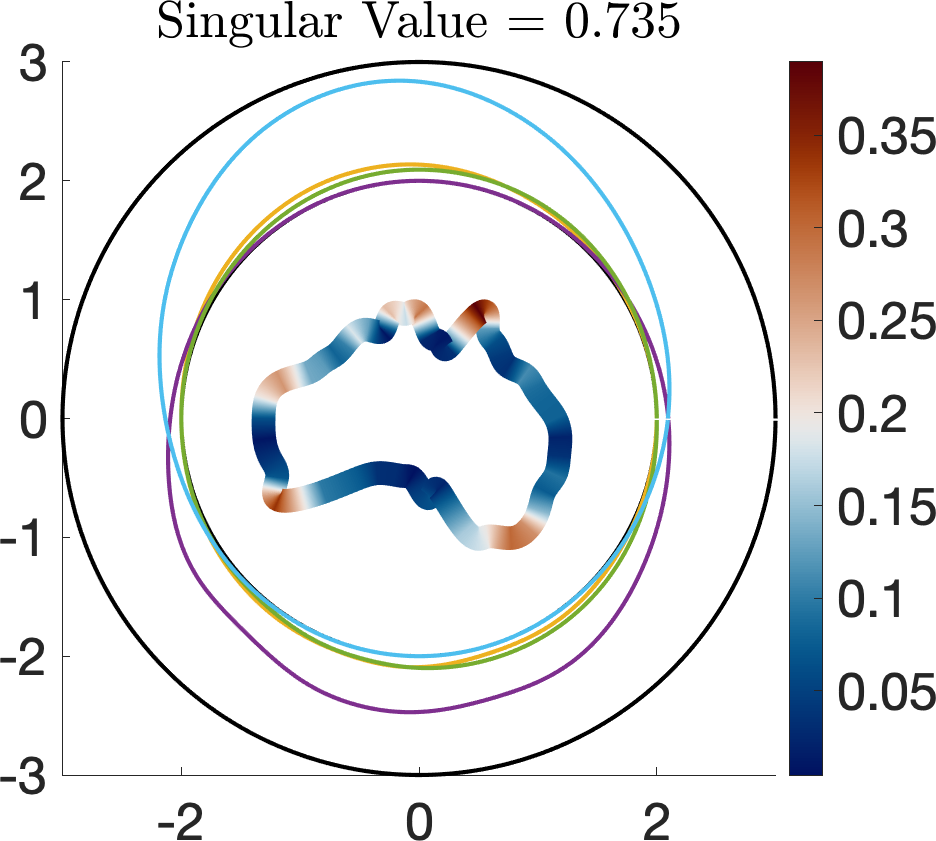}}; \\
   };
\end{tikzpicture}
\caption{Visualisation of the first 3 singular values and associated right and left singular vectors for $N_\text{spline}=48$ and $|\mathbf{J}|$ comprising of $4$ incident waves with $\mathbf{\widehat{{d}}}= \pm \boldsymbol{\hat{e}_x}, \pm \boldsymbol{\hat{e}_y}$.}
\label{fig:4 inc wave SVD N=48}
\end{figure}

Next we consider the case of $4$ incident waves with $\mathbf{\widehat{{d}}}= \pm \boldsymbol{\hat{e}_x}, \pm \boldsymbol{\hat{e}_y}$. 
In Figures \ref{fig:4 inc wave SVD N=12}, \ref{fig:4 inc wave SVD N=48}
we visualise the left (shape) and right (far-field) singular vectors
for the largest three singular values
for $ka=2\pi, 4\pi, 6\pi$ respectively.
We note that increasing $ka$ dictates a slower rate of decay of the singular values.
We visualise the shape singular vectors $|\mathbf{u_i}|$ 
using colour on the scatterer boundary.
For the far-field singular vectors $|\mathbf{v_i}|$ visualisation, we plot the associated acoustic cross section $|\mathbf{v_i}|^2$ 
in polar coordinates. Since we have 4 incident waves, there are 4 associated cross sections, discriminated
by color.
The far-field singular vectors are
normalised by $\sigma_i$ so that the size of the perturbation in the far field
for a given shape perturbation can be quantified.
For example, in the top right plots of both figures the far field is close to zero
(the inner circle) indicating that the far field is not very sensitive to the
third singular shape-vector.

These modes provide insight into which boundary features are most visible from far-field data across all incident waves. 
The most noticeable feature is that the support of the first singular shape-vector 
is in all cases associated with the 
south-western tip of Western Australia. The corresponding singular far-field vectors
show that this feature produces large perturbations in the far field for those incident
waves that illuminate it.

All rows of Figure \ref{fig:4 inc wave SVD N=12} show considerable sensitivity across the boundary except for the concave region in the southern tip corresponding to $\theta_{10}^{\partial D}$. In turn, this means that under these conditions, two scatterers
that are identical  except for small deviations in $\omega_{10}$ may have indistinguishable far-fields patterns, manifesting ill-posed behaviour. In Figure \ref{fig:4 inc wave SVD N=48} the situation is significantly more complex, with many more identifiable regions with small sensitivity. This is nevertheless expected, since the finer details in the scatterer would require much higher $ka$ values due to the \textit{diffraction limit}.

\section{Conclusions} \label{section:conclusions}
We have provided an investigation of the effects of local perturbations to a given $2$D scatterers shape in corresponding far-field patterns. These effects can be written in terms of the Fr\'echet derivative of the far-field operator $\mathrm{d} \mathcal{F}[\cdot]/\mathrm{d} r$ with $N_\text{spline}$ localised input functions following from a splines representation for the scatterer. The computation reduces to a coupled BIE which is accurately solved using the Nystr\"om method. Upon choosing observation directions, all results can be represented through a sensitivity matrix $\mathbf{J}$ which we show can conveniently be studied using its SVD. This tool gives us instant useful information 
into the ill-posedness of the inverse problem
such as which shape deformations can be reliably recovered from far-field data and, on the contrary, which are effectively invisible to far-field measurements.

\bibliographystyle{unsrt}

\end{document}